\documentclass[final,a4paper,onecolumn,3p]{elsarticle}
\pdfoutput=1
\usepackage[colorlinks=true,linkcolor=blue,citecolor=blue,urlcolor=blue]{hyperref}
\usepackage{lineno}
\usepackage[utf8]{inputenc}
\usepackage[english]{babel}
\usepackage{amsmath}
\usepackage{amsfonts}
\usepackage{amssymb}
\usepackage{graphicx}
\usepackage{epstopdf}
\usepackage{listings}
\usepackage{empheq}
\usepackage{url}

\usepackage{listings}
\usepackage{xcolor}
\definecolor{codegreen}{rgb}{0,0.6,0}
\definecolor{codegray}{rgb}{0.5,0.5,0.5}
\definecolor{codepurple}{rgb}{0.58,0,0.82}
\definecolor{backcolour}{rgb}{0.95,0.95,0.92}
\lstdefinestyle{mystyle}{
	backgroundcolor=\color{backcolour},   
	commentstyle=\color{codegreen},
	keywordstyle=\color{blue},
	numberstyle=\tiny\color{codegray},
	stringstyle=\color{brown},
	basicstyle=\ttfamily\footnotesize,
	breakatwhitespace=false,         
	breaklines=true,                 
	captionpos=b,                    
	keepspaces=true,                 
	numbers=left,    
	numbersep=5pt,                  
	showspaces=false,                
	showstringspaces=false,
	showtabs=false,                  
	tabsize=2
}
\usepackage{fancyhdr}
\begin{document}
\renewcommand{\today}{26 January 2020}

\begin{abstract}
We present a novel application of the recently developed AAA algorithm to the solution of Laplace 2D problems; an application to conformal mapping is also shown as a particular case. These classes of problems have also been addressed by means of dedicated software implementations based on rational function approximation, with unprecedented speed and accuracy. Still we show how the AAA algorithm, conceived to be a flexible and domain independent method, is capable of addressing these situations within a unified framework.
\end{abstract}


\begin{frontmatter}
	\title{Solving Laplace problems with the AAA algorithm}
	
	\author{Stefano Costa\fnref{fn1}}
	\ead{stefano.costa@ieee.org}
	
	\fntext[fn1]{IEEE, member; Piacenza 29122 Italy. {\includegraphics[scale=0.5,keepaspectratio]{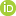}} \url{ https://orcid.org/0000-0002-9434-356X}}
	
	\journal{arXiv}
\end{frontmatter}

\section{Introduction}
Recently, problems of the Laplace type and conformal mapping, which can be shown to reduce to a particular form of a Laplace problem \cite{Tr19}, have been addressed by the strikingly efficient software implementation based on rational approximations called the "lightning Laplace solver" (LLS) in domain with corners \cite{GT19b}, and the "Vandermonde with Arnoldi" algorithm for the stabilization of matrix computation associated with power sums in regions with smooth boundaries \cite{Tr19b}. The relevant literature also shows how it is crucial that when singularities like corners are present, sample points be exponentially clustered in their neighbourhood to achieve root-exponential convergence \cite{GT19a,GT19b}. Our method starts from the aim of giving a solution also for those cases, frequently found in physics and engineering, where only a set of values associated to a generic "reasonably discretized" contour is given, and results with an accuracy of a fraction of a percent are sufficient. Being based on the flexible Adaptive Antoulas-Anderson (AAA) algorithm \cite{NOT19}, it can help whenever parametric boundaries or explicit function definitions are not available.

This paper is focused on coding. After a brief presentation of the method, illustrating the underlying idea and its mathematical definition, we give a collection of selected problems along with the listings for their solution.

\section{Method and Algorithm}
\begin{figure}[!h]
	\centering
	\includegraphics[width=1\linewidth]{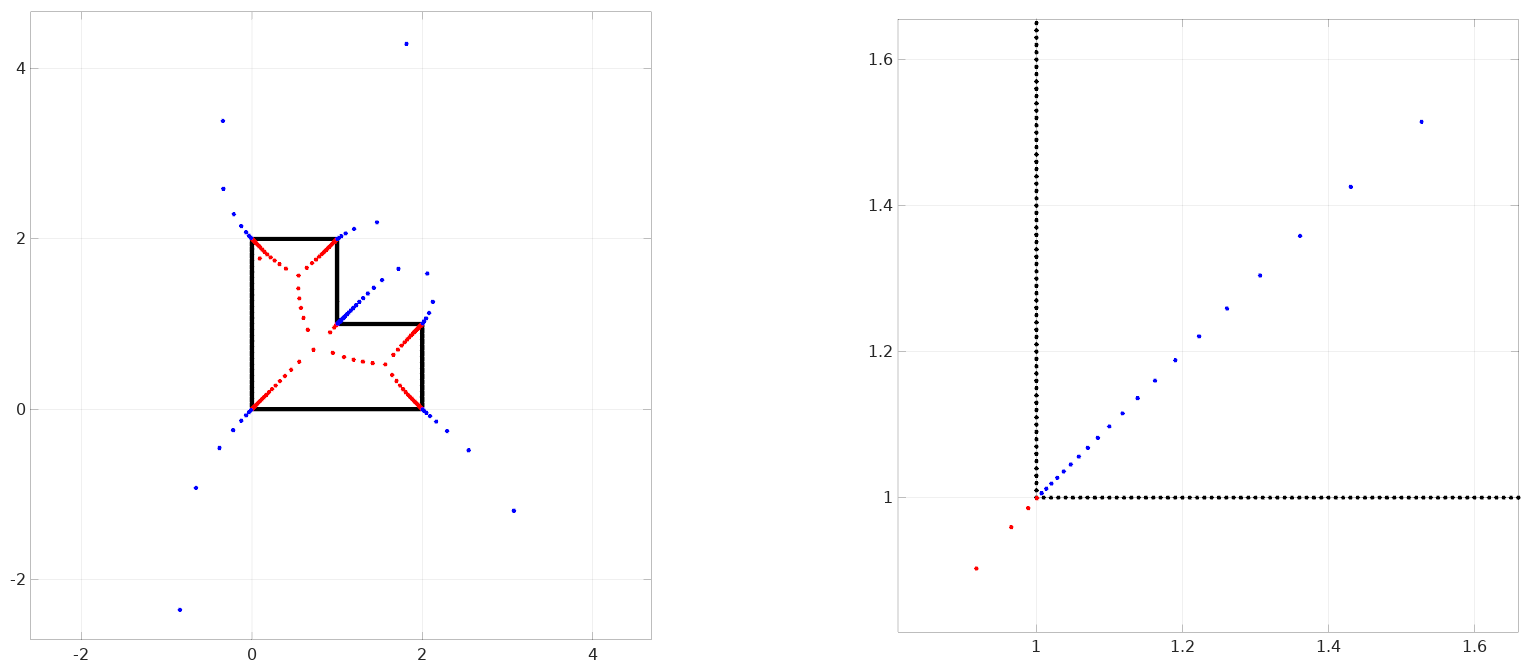}
	\caption{Poles returned by the AAA approximation  of $u(z)=(\Re z)^2=x^2$ on the L-shaped boundary (left) and closeup of the reentrant corner (right); note the pole clustering notwithstanding the even boundary sampling.}
	\label{fig:lpoles}
\end{figure}

The AAA algorithm computes rational barycentric approximations of scalar functions defined on subsets of the complex plane, also returning their poles, residues and zeros: this allows representations in terms of partial fractions. Figure \ref{fig:lpoles} shows the poles returned after the approximation of $u(z)=(\Re z)^2=x^2$ on an L-shaped boundary discretized by evenly sampled points.
We first notice that poles are arranged in a somewhat exponentially clustered fashion on both sides of the boundary to represent corners, irrespective of sampling frequency; the trouble is that they are spread everywhere with no chance to exclude them from a particular region. Roughly speaking, if the whole set provides a tolerance of e.g. 1e-13, half their number should achieve about 1e-6, which is still decent indeed, so the idea is as simple as this: if we want to solve a Laplace problem with Dirichlet boundary condition $u(z)$ in the inside region, we must drop the inner poles; if we want to solve outside, then we must drop the outer ones: in other words, we drop the subset impairing analyticity. Clearly this plain concept alone is not sufficient, as trivial examples can show, therefore a replacement for the excluded subset need be found.

Our method, from a formal viewpoint, has much in common with the LLS fully described in \cite{GT19b}: given a real-valued function $u(z)$ defined on a contour $\gamma$ in the complex plane, and dividing it into an interior region $\Gamma_\oplus$ and an exterior one $\Gamma_\ominus$, find a function $w(z)$ analytic in $\Gamma_\oplus$
\begin{align}
	\begin{matrix}
		w(z) = & \underbrace{\sum_{n=0}^{N} a_n(z-c)^n} & + & \underbrace{b_0 + \sum_{p_k\in\Gamma_\ominus}\dfrac{b_k}{z-p_k}} \\
		& smooth,\;\tilde{w}(z) & & singular,\;\hat{w}(z)
	\end{matrix}
	\label{eqn:main}
\end{align}
for which $\Re w(z)\approx u(z)$ and $\Im w(z)\approx v(z)$: often in applications these are referred to as the \emph{scalar potential} and the \emph{stream function} respectively, of the \emph{complex potential} $w(z)$. $c$ is the centre of a power sum describing the "smooth" part of the problem, placed in about the middle of the geometry, whilst $\{p_k\}$ is the subset of poles in $\Gamma_\ominus$ shaping the "singular" part, namely contour corners and irregular function behaviours. The interior region is defined with respect to $\gamma$: we let $\Gamma_\oplus$ be the region located on the left of the boundary when traversing it counterclockwise. Conversely, if we want a solution in $\Gamma_\ominus$ we take $\{p_k\}\in\Gamma_\oplus$ and $N$ with a \emph{negative} value in formula (\ref{eqn:main}), $c$ resulting in an additional arbitrary pole placed in the interior region. The magnitude of $N$ is discussed later.

\begin{figure}[!h]
	\centering
	\includegraphics[width=1\linewidth]{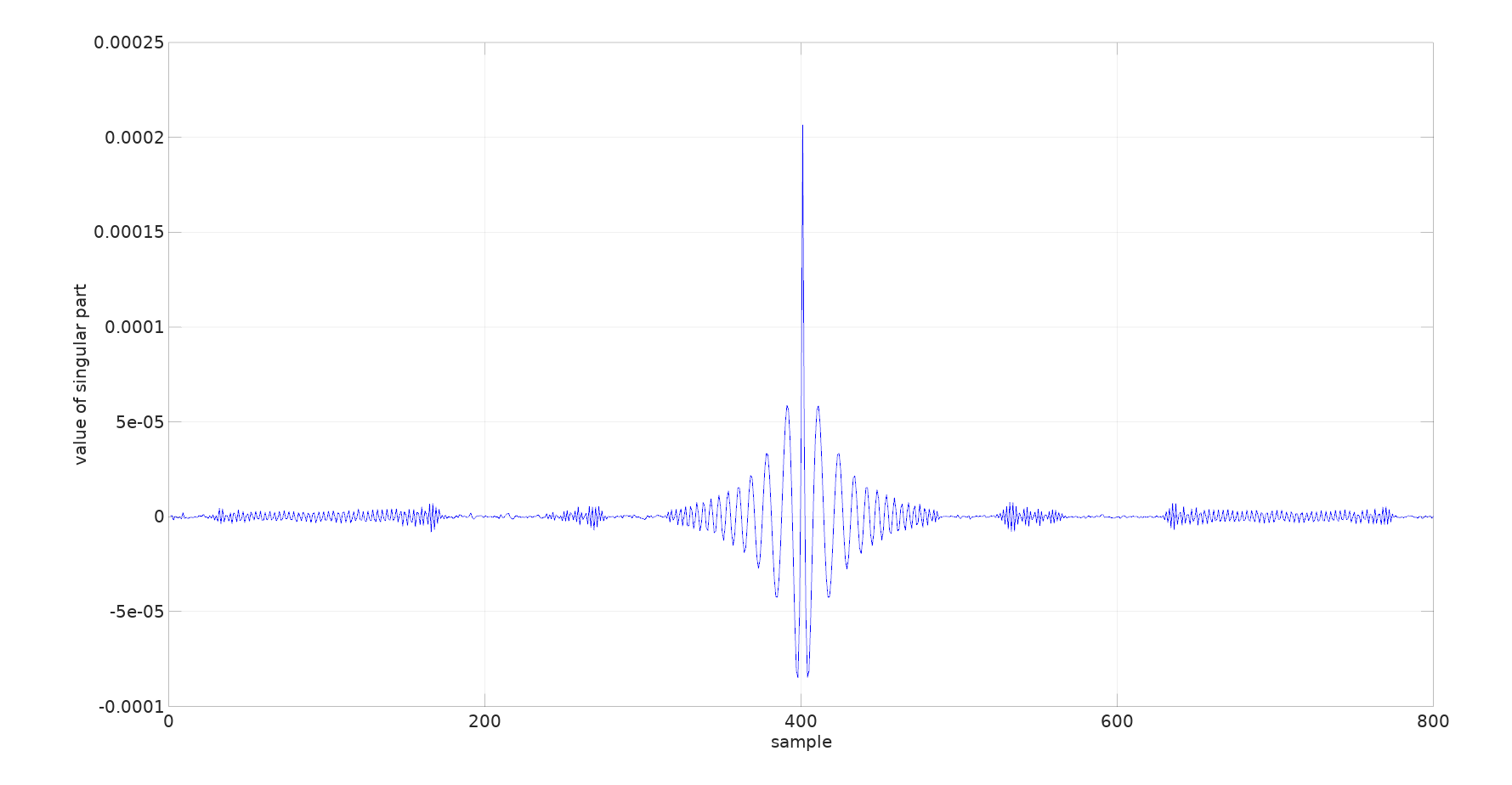}
	\caption{Choosing a large value for $N$ in (\ref{eqn:main}), far from being beneficial, typically causes $\hat{u}(z)$ to behave as shown: widely oscillating and "noisy". In such situations the AAA algorithm is unable to return sensible approximations, even if a very large degree is allowed.}
	\label{fig:ntoolarge}
\end{figure}

From now on, we pursue a different path from \cite{GT19b}: the idea is to build the two parts $\tilde{w}(z)$ and $\hat{w}(z)$ in sequence. First, a point $c$ is chosen which lies about in the middle of the region, and $\{a_n=\alpha_n+i\beta_n\}$ are determined by solving the least-squares (LS) problem
\begin{align}
	\mathbf{A}_r = \begin{bmatrix}
		A_{j,n}=\Re(z_j-c)^n
	\end{bmatrix}, \\
	\mathbf{A}_i = \begin{bmatrix}
		A_{j,n}=\Im(z_j-c)^n
	\end{bmatrix}, \\
	\begin{bmatrix}
		\mathbf{A}_r & -\mathbf{A}_i
	\end{bmatrix}
	\begin{bmatrix}
		\alpha_n \\
		\beta_n
	\end{bmatrix} \approx u(z_j)
\end{align}
where $\{z_j\}$ is the given set of points describing $\gamma$, and $\{u(z_j)\}$ the corresponding function values. The real singular part  $\hat{u}(z)$ of the problem is pulled out by simply subtracting
\begin{align}
	\hat{u}(z_j) = u(z_j)-\Re\sum_{n=0}^{N} a_n(z_j-c)^n = u(z_j)-\tilde{u}(z_j).
\end{align}
Then $\hat{u}(z)$ is approximated by the AAA algorithm, unwanted poles are dropped,  and a second LS problem is solved in order to determine $\{b_k=\delta_k+i\epsilon_k\}$:
\begin{align}
	\mathbf{B}_r = \begin{bmatrix}
			B_{j,k}=\Re(1/(z_j-p_k))
		\end{bmatrix}, \\
	\mathbf{B}_i = \begin{bmatrix}
			B_{j,k}=\Im(1/(z_j-p_k))
		\end{bmatrix}, \\
		\begin{bmatrix}
			\mathbf{B}_r & -\mathbf{B}_i & \mathbf{1}
		\end{bmatrix}
		\begin{bmatrix}
			\delta_k \\
			\epsilon_k \\
			b_0
		\end{bmatrix} \approx \hat{u}(z_j).
\end{align}
The reason for there being two constants in (\ref{eqn:main}) is that $b_0$ adjusts the value of $a_0$ in the smooth part, in order to approximate $\{u(z_j)\}$ better with the second LS fitting. The parameter $N$ in (\ref{eqn:main}) deserves particular attention. One usually works with sample points (almost) regularly spaced in applications, but picking a value with the same order of magnitude as sample number would imply an overwhelming amount of terms with unacceptable computing time. Even worse, too large a number would compromise the subsequent approximation of the singular part: the difference between $u(z)$ and its first LS fitting $\tilde{u}(z)$ along the boundary would be a function severely affected by the Gibbs phenomenon, not restricted to corner neighbourhoods. Figure \ref{fig:ntoolarge} illustrates such situation, where $\hat{u}(z)$ can not be approximated correctly by the AAA algorithm. Therefore, $N$ is determined by taking the logarithm of sample number: it is an engineering choice, typically giving rise to some 10 terms in the smooth part, and we don’t claim it to be optimal though in a number of cases it has appeared to work well. In cases of concave shapes close to a circle (a square, a regular polygon, an oval, or the like) an even smaller value need be used.

We are now ready to give an overview of the algorithm. We exploit the recently announced Arnoldi stabilization for Vandermonde matrices \cite{Tr19b}, which provides superior stability to the ill-conditioned computation of the smooth part in the first step:
\begin{enumerate}
	\item Approximate the smooth part $\tilde{u}(z)$ by means of a power sum, analytic in the desired region;
	\item Pull the singular part $\hat{u}(z)$ out of $u(z)$, and approximate it by means of the AAA algorithm;
	\item Use the poles returned by the AAA approximation in the other region to solve the second LS problem for $\hat{u}(z)$;
	\item Build a function handle for $w(z)$, and estimate the maximum error $|u-\Re w|_\infty$ on the boundary.
\end{enumerate}
In the next sections, detailed solutions to some problems of the Laplace type and conformal mapping are given. Some general remarks:
\renewcommand{\labelitemi}{$\square$}
\begin{itemize}
	\item Large portions of the code are repeated in different listings: we always give the full solution ready for cut-and-paste into the computing environment. This widespread uniformity and standardization, we believe, will make the reader appreciate the generality of the method.
	\item Boundary is always represented starting from a uniform discretization step, to give evidence of the AAA algorithm's capability of clustering poles as necessary, thus partially compensating any lack of exponential sampling close to critical points.
	\item The AAA algorithm is called with a high maximum degree of 1000 as the default value of 100 is seldom sufficient to reach the default tolerance of 1e-13. In all the proposed examples an approximation with about 200 poles is obtained;
	\item Lawson iterations for minimax approximation are very likely to cause poles to reposition badly, so they must be turned off. If doubts arise about the possibility of getting a number of zeros greater than the number of poles, one iteration only should be used.
	\item The maximum error $|u-\Re w|_\infty$ on the boundary is computed for available samples, and accuracy can not be confirmed by checking on a twice or four times finer mesh as the LLS does, therefore in general it can not be considered an accurate upper bound. Experiments show that a much bigger error is to be expected throughout the domain if no particular care is taken in the neighbourhood of critical points.
\end{itemize}
The following are needed to run the examples:
\begin{itemize}
	\item A MATLAB 6 or Octave 4.2 environment, or newer;
	\item The Chebfun package \cite{DHT14}, carrying the most recent version of the AAA algorithm;
	\item The \texttt{polyfitA.m} and \texttt{polyvalA.m} functions described in \cite{Tr19b}, for stable computations of power sums; these must be found in the environment path;
	\item The SC-Toolbox package \cite{SCTB} for the \texttt{polygon()} and \texttt{isinpoly()} functions, to build a polygon and check whether points fall within it.
\end{itemize}

\newpage
\section{Laplace problems with Dirichlet boundary condition $u(z)=(\Re z)^2$}

\subsection{Analytic solution in the L-shaped domain interior}
\begin{lstlisting}[language=Matlab, caption={}]
stp = 0.01;                                             % boundary step discretization
c = 0.5+0.5*i;                                          % power sum centre

z = [0:stp:2]; z = [z 2+i*[0:stp:1]];                   % L-shaped boundary definition
z = [z [2:-stp:1]+i]; z = [z 1+i*[1:stp:2]];
z = [z [1:-stp:0]+2*i]; z = [z i*[2:-stp:0]].';
p = polygon([0 2 2+i 1+i 1+2*i 2*i]);

u = real(z).^2;                                         % Dirichlet boundary condition

N = 10+ceil(log(length(z)));                            % smooth part, nr. of terms
[a,H] = polyfitA(z-c,u,N);                              % Vandermonde with Arnoldi
wz = @(x) polyvalA(a,H,x(:)-c);                         % function handle, smooth part
kz = imag(wz(c));

up = u-real(wz(z));                                     % singular part of b.c.
[r,pol] = aaa(up,z,'mmax',1000,'lawson',0);             % AAA approximation
m = find(isinpoly(pol,p,1E-16)==0);                     % find poles outside
polm = pol(m);
B = 1./(z-polm.');                                      % coefficient matrix
B(:,end+1) = 1;
B = [real(B) -imag(B)];
b = reshape(B\up,[],2);                                 % LS approximation
b = b(:,1)+i*b(:,2);
wp = @(x) [1./(x(:)-polm.') ones(size(x(:)))]*b;        % function handle, singular part
kp = imag(wp(c));

w = @(x) wp(x)+wz(x(:))-i*(kz+kp);                      % function handle, solution
maxerr = norm(u-real(w(z)),'inf')                       % maximum error at the boundary
\end{lstlisting}
\begin{figure}[!h]
	\centering
	\includegraphics[width=1\linewidth]{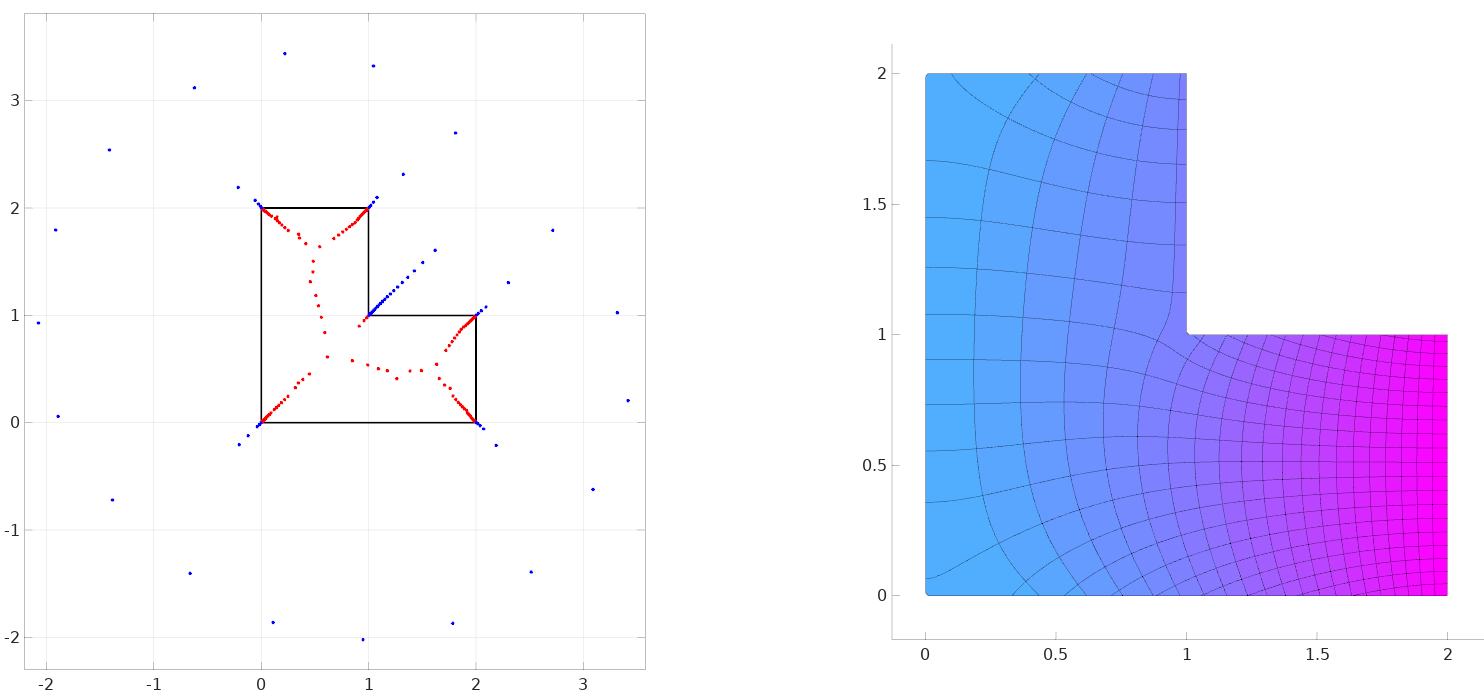}
	%\caption{}
	\label{fig:lap1_l}
\end{figure}
Field solution with uniform discretization step of 0.01 at the boundary, $N$ = 17 in the smooth part and 186 poles, of which 66 fall in the exterior region. $|u-\Re w|_\infty$ = 6.47e-7 close to that returned by the LLS when called with default parameters, even though this can not be considered a true upper bound for the error.

\begin{figure}[!h]
	\centering
	\includegraphics[width=1\linewidth]{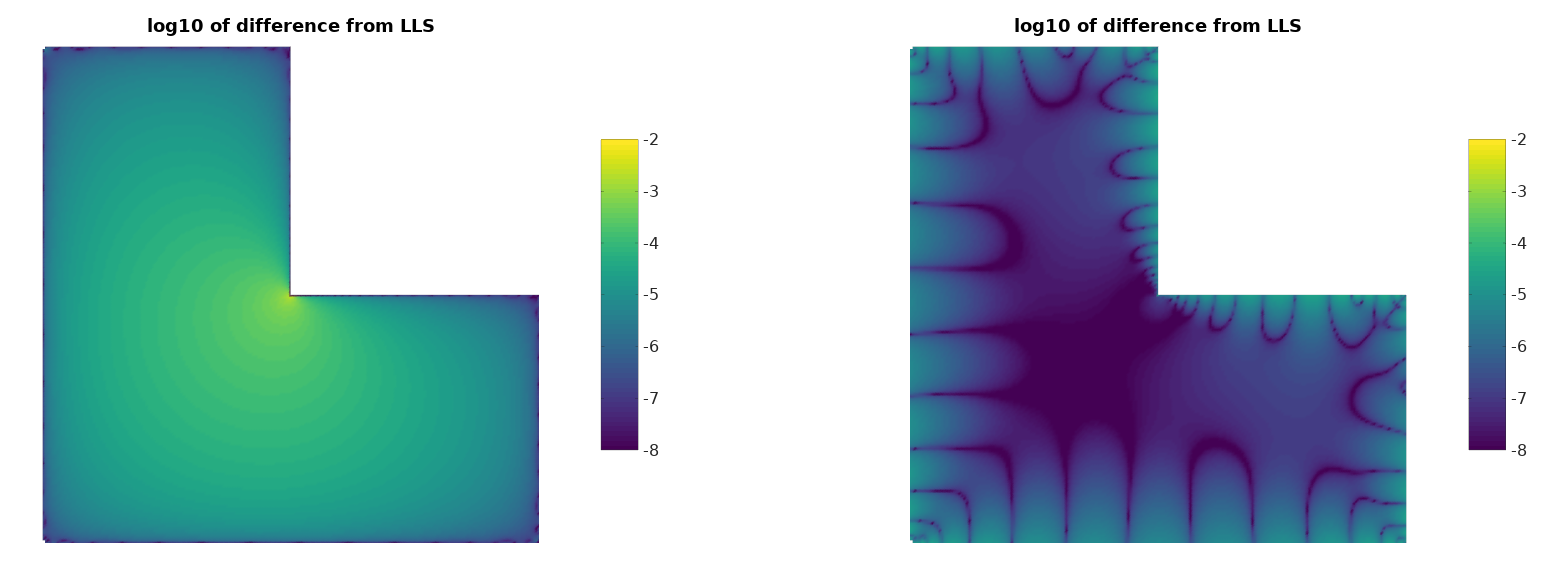}
	%\caption{}
	\label{fig:laperr_l}
\end{figure}
A comparison throughout the domain with the values of $u(z)$ computed by the LLS, called with a tolerance for maximal absolute error of 1e-13, and therefore considered "exact", shows clearly on the left a possible drawback of picking evenly sampled points along the boundary. In this case an excellent pointwise boundary approximation does not override the high local error at the reentrant corner, even though in this very point $|u-\Re w|_\infty$ is very low. As a testbed \cite{Tr18} consider the computation of $u(0.99+0.99\,i)$, the exact value being 1.0267919261073... We get $\Re w(0.99+0.99\,i)$ = 1.025 (-0.16\%). However, consider that, as reported in \cite{GT19c}, experts using specialized FEM codes gave solutions to 2-4 digits of accuracy, and it is still much better than values obtained from many commercial and open source FEM packages, usually affected by errors between 1\% and 10\% if no particular attention is paid to mesh generation near critical points.

On the right, the solution after adding 50 clustered points per side for each corner (i.e. 100 points/corner) to the evenly spaced set. This time we get $\Re w(0.99+0.99\,i)$ = 1.0267918, and $|u-\Re w|_\infty$ = 4.82e-5 represents a reliable upper bound. Additional points are obtained from interpolation of the given samples very easily, by means of calling \texttt{logspace(-6,0,50)} to determine their displacement from corners along the boundary. In other words they start from a distance of 1e-6 from corners and proceed logarithmically spaced. The degree of the AAA approximation is limited to 200 to gain in speed with no drawbacks, as in this case clustered points make the difference. Whenever precise interpolations can be computed from given sets of values, this simple technique is recommended to improve final results greatly.

\newpage
\subsection{Analytic solution in the blade-shaped domain interior,  $\gamma=2\cos\vartheta+i(\sin\vartheta+2\cos^3\vartheta)$}
\begin{lstlisting}[language=Matlab, caption={}]
h = linspace(0,2*pi,500)';                              % discretization, angle theta
c = 0;                                                  % power sum centre
z = 2*cos(h)+i*(sin(h)+2*cos(h).^3);                    % blade-shaped boundary definition
p = polygon(z);

u = real(z).^2;                                         % Dirichlet boundary condition

N = 10+ceil(log(length(z)));                            % smooth part, nr. of terms
[a,H] = polyfitA(z-c,u,N);                              % Vandermonde with Arnoldi
wz = @(x) polyvalA(a,H,x(:)-c);                         % function handle, smooth part
kz = imag(wz(c));

up = u-real(wz(z));                                     % singular part of b.c.
[r,pol] = aaa(up,z,'mmax',1000,'lawson',0);             % AAA approximation
m = find(isinpoly(pol,p,1E-16)==0);                     % find poles outside
polm = pol(m);
B = 1./(z-polm.');                                      % coefficient matrix
B(:,end+1) = 1;
B = [real(B) -imag(B)];
b = reshape(B\up,[],2);                                 % LS approximation
b = b(:,1)+i*b(:,2);
wp = @(x) [1./(x(:)-polm.') ones(size(x(:)))]*b;        % function handle, singular part
kp = imag(wp(c));

w = @(x) wp(x)+wz(x(:))-i*(kz+kp);                      % function handle, solution
maxerr = norm(u-real(w(z)),'inf')                       % maximum error at the boundary
\end{lstlisting}
\begin{figure}[!h]
	\centering
	\includegraphics[width=1\linewidth]{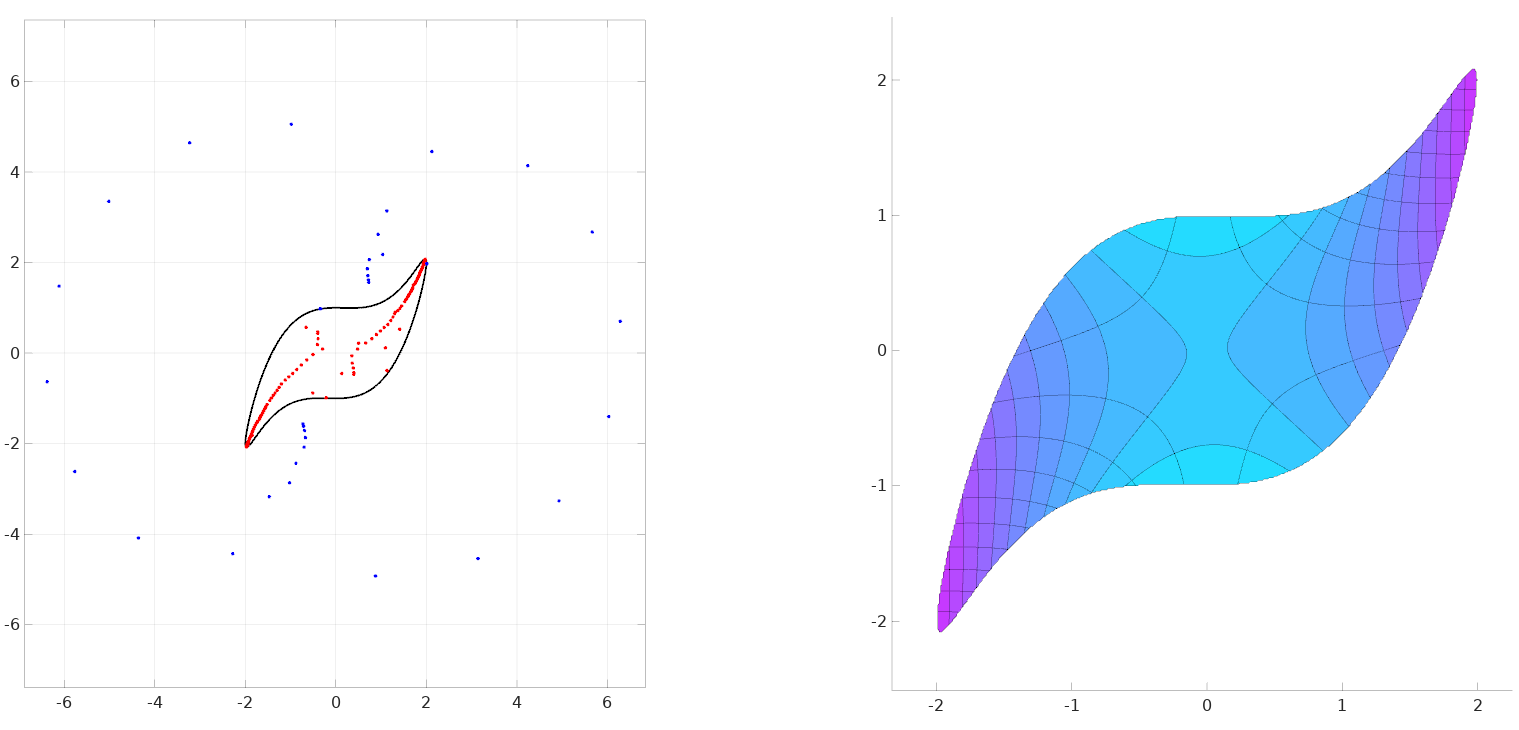}
	%\caption{}
	\label{fig:lap2_bl}
\end{figure}
Field solution with $N$ = 17 in the smooth part and 167 poles, of which 34 fall in the exterior region. The smooth contour is represented by a polygonal union of tiny segments, getting smaller at the blade ends. $|u-\Re w|_\infty$ = 4.10e-7. The listing is the same as in the previous case except in the boundary definition.

\newpage
\subsection{Analytic solution in the L-shaped domain exterior}
\begin{lstlisting}[language=Matlab, caption={}]
stp = 0.01;                                             % boundary step discretization
c = 0.5+0.5*i;                                          % power sum centre

z = [0:stp:2]; z = [z 2+i*[0:stp:1]];                   % L-shaped boundary definition
z = [z [2:-stp:1]+i]; z = [z 1+i*[1:stp:2]];
z = [z [1:-stp:0]+2*i]; z = [z i*[2:-stp:0]].';
p = polygon([0 2 2+i 1+i 1+2*i 2*i]);

u = real(z).^2;                                         % Dirichlet boundary condition

N = 10+ceil(log(length(z)));                            % smooth part, nr. of terms
[a,H]=polyfitA(1./(z-c),u,N);                           % Vandermonde with Arnoldi
wz = @(x) polyvalA(a,H,1./(x(:)-c));                    % function handle, smooth part
kz = imag(wz(inf));

up = u-real(wz(z));                                     % singular part of b.c.
[r,pol] = aaa(up,z,'mmax',1000,'lawson',0);             % AAA approximation
m = find(isinpoly(pol,p,1E-16)==1);                     % find poles inside
polm = pol(m);
B = 1./(z-polm.');                                      % coefficient matrix
B(:,end+1) = 1;
B = [real(B) -imag(B)];
b = reshape(B\up,[],2);                                 % LS approximation
b = b(:,1)+i*b(:,2);
wp = @(x) [1./(x(:)-polm.') ones(size(x(:)))]*b;        % function handle, singular part
kp = imag(wp(inf));

w = @(x) wp(x)+wz(x(:))-i*(kz+kp);                      % function handle, solution
maxerr = norm(u-real(w(z)),'inf')                       % maximum error at the boundary
\end{lstlisting}
\begin{figure}[!h]
	\centering
	\includegraphics[width=1\linewidth]{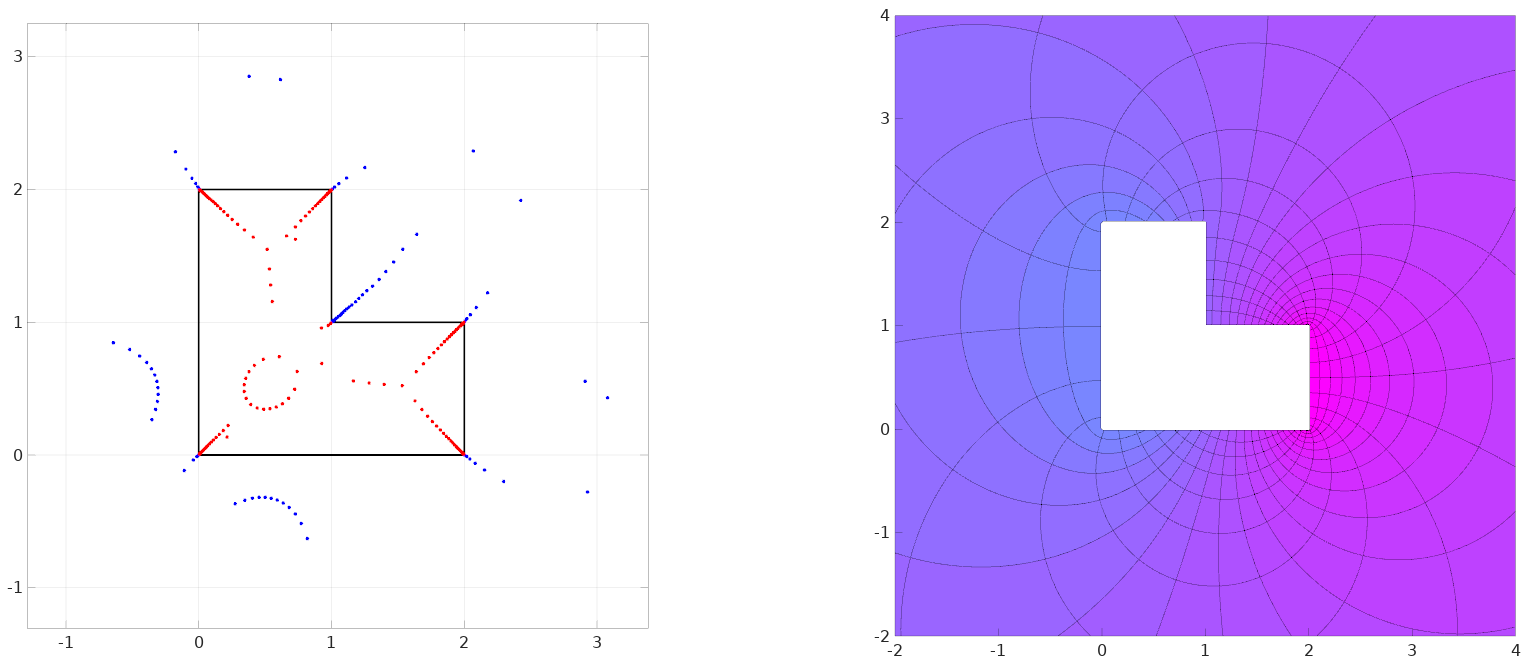}
	%\caption{}
	\label{fig:lap3_l}
\end{figure}
Field solution with $N$ = 17 in the smooth part, meaning $N$ = -17 in the first power sum of (\ref{eqn:main}), and 205 poles, of which 124 fall in the interior region. $|u-\Re w|_\infty$ = 3.07e-5. The listing differs from the first case (solution in the interior region) in the definition of the smooth part, and in the \texttt{isinpoly()} function asking for interior poles.

\section{Conformal mapping}

\subsection{L-shaped domain interior onto/from disk interior}
\begin{lstlisting}[language=Matlab, caption={}]
stp = 0.01;                                             % boundary step discretization
c = 0.5+0.5*i;                                          % power sum centre

z = [0:stp:2]; z = [z 2+i*[0:stp:1]];                   % L-shaped boundary definition
z = [z [2:-stp:1]+i]; z = [z 1+i*[1:stp:2]];
z = [z [1:-stp:0]+2*i]; z = [z i*[2:-stp:0]].';
p = polygon([0 2 2+i 1+i 1+2*i 2*i]);

u = -log(abs(z-c));                                     % Dirichlet boundary condition

N = 10+ceil(log(length(z)));                            % smooth part, nr. of terms
[a,H] = polyfitA(z-c,u,N);                              % Vandermonde with Arnoldi
wz = polyvalA(a,H,z-c);                                 % smooth part

up = u-real(wz);                                        % singular part of b.c.
[r,pol] = aaa(up,z,'mmax',1000,'lawson',0);             % AAA approximation
m = find(isinpoly(pol,p,1E-16)==0);                     % find poles outside
B = 1./(z-pol(m).');                                    % coefficient matrix
B(:,end+1) = 1;
V = [real(B) -imag(B)];
b = reshape(V\up,[],2);                                 % LS approximation
b = b(:,1)+i*b(:,2);
wp = B*b;                                               % singular part

w = wp+wz;                                              % solution
maxerr = norm(u-real(w),'inf')                          % maximum error at the boundary

g = (z-c).*exp(w);                                      % conformal mapping function
f = aaa(g,z);                                           % polygon onto disk, AAA approx.
finv = aaa(z,g);                                        % disk onto polygon, AAA approx.
\end{lstlisting}
\begin{figure}[!h]
	\centering
	\includegraphics[width=1\linewidth]{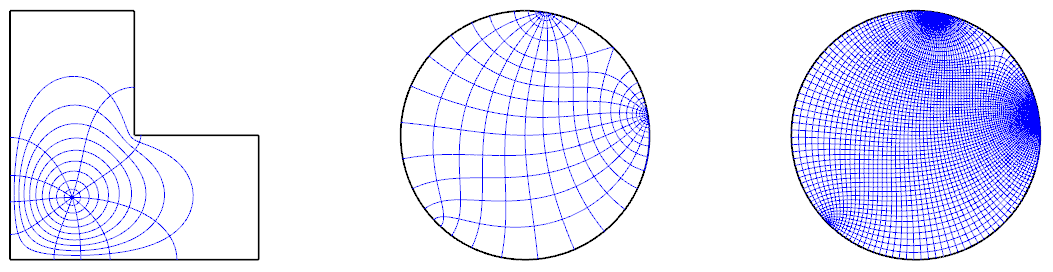}
	%\caption{}
	\label{fig:cm1_l2}
\end{figure}
Conformal mapping (see \cite{Tr19,Tr19b}) with uniform discretization step of 0.01 at the boundary: from disk onto polygon (mapping of concentric circles, left) and from polygon onto disk (mapping of square grids of different size, centre and right); this  $N$ = 17 in the smooth part and 192 poles, of which 71 fall in the exterior region. $|u-\Re w|_\infty$ = 5.62e-8. Even though sample points are not exponentially clustered near corners, an acceptable accuracy is maintained: images of grid lines at a distance of 0.01 from the L-shaped boundary match those computed with the SC-Toolbox within some 1e-3.

\newpage
\subsection{L-shaped domain exterior onto/from disk interior}
\begin{lstlisting}[language=Matlab, caption={}]
stp = 0.01;                                             % boundary step discretization
c = 0.5+0.5*i;                                          % power sum centre

z = [0:stp:2]; z = [z 2+i*[0:stp:1]];                   % L-shaped boundary definition
z = [z [2:-stp:1]+i]; z = [z 1+i*[1:stp:2]];
z = [z [1:-stp:0]+2*i]; z = [z i*[2:-stp:0]].';
p = polygon([0 2 2+i 1+i 1+2*i 2*i]);

u = log(abs(z-c));                                      % Dirichlet boundary condition

N = 10+ceil(log(length(z)));                            % smooth part, nr. of terms
[a,H] = polyfitA(1./(z-c),u,N);                         % Vandermonde with Arnoldi
wz = polyvalA(a,H,1./(z-c));                            % smooth part

up = u-real(wz);                                        % singular part of b.c.
[r,pol] = aaa(up,z,'mmax',1000,'lawson',0);             % AAA approximation
m = find(isinpoly(pol,p,1E-16)==1);                     % find poles inside
B = 1./(z-pol(m).');                                    % coefficient matrix
B(:,end+1) = 1;
V = [real(B) -imag(B)];
b = reshape(V\up,[],2);                                 % LS approximation
b = b(:,1)+i*b(:,2);
wp = B*b;                                               % singular part

w = wp+wz;                                              % solution
maxerr = norm(u-real(w),'inf')                          % maximum error at the boundary

g = (1./(z-c)).*exp(w);                                 % conformal mapping function
f = aaa(g,z);                                           % polygon onto disk, AAA approx.
finv = aaa(z,g);                                        % disk onto polygon, AAA approx.
\end{lstlisting}
\begin{figure}[!h]
	\centering
	\includegraphics[width=0.6\linewidth]{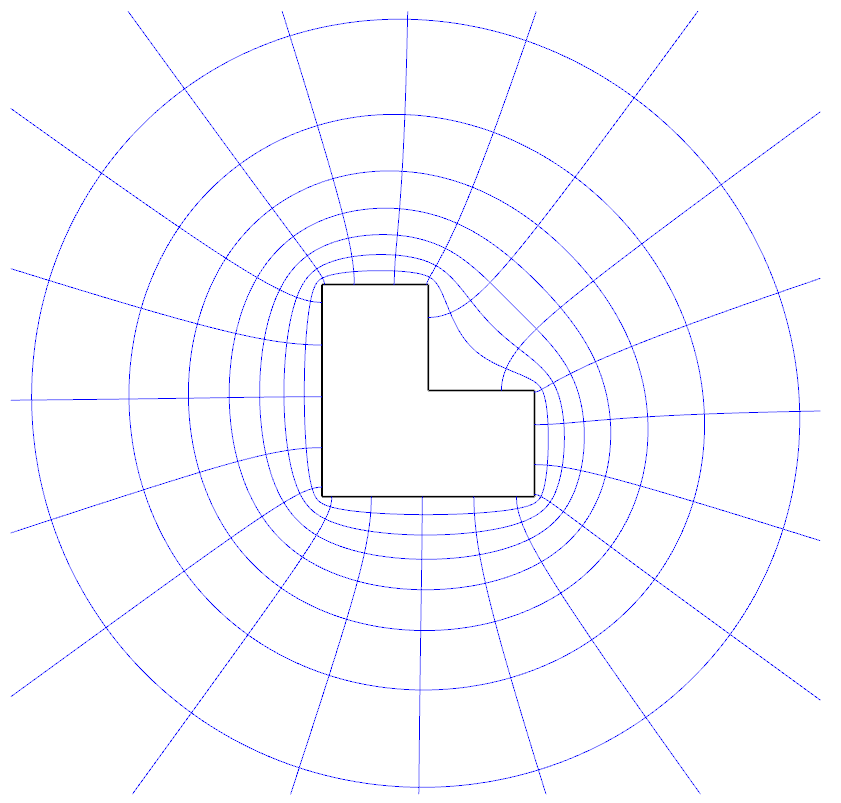}
	%\caption{}
	\label{fig:cm2_l}
\end{figure}
Conformal mapping with $N$ = 17 in the smooth part and 204 poles, of which 124 fall in the interior region. $|u-\Re w|_\infty$ = 7.61e-7.

\newpage
\subsection{Doubly connected domain onto/from annulus}
\begin{lstlisting}[language=Matlab, caption={}]
stp = 0.01;                                             % boundary step discretization
c = -0.25-0.25*i;                                       % power sum centre;

s1 = [0:stp:2];                                         % outer Jordan curve
z1 = [i*sqrt(2)+s1*exp(i*5/4*pi)]; z1 = [z1 -sqrt(2)+s1*exp(i*7/4*pi)];
z1 = [z1 -i*sqrt(2)+s1*exp(i*1/4*pi)]; z1 = [z1 sqrt(2)+s1*exp(i*3/4*pi)].';
p1 = polygon([i*sqrt(2) -sqrt(2) -i*sqrt(2) sqrt(2)]);

s2 = [0:stp:0.5];                                       % inner Jordan curve
z2 = s2-0.5-0.5*i; z2 = [z2 s2*exp(i*pi/2)+0.0-0.5*i];
z2 = [z2 s2*exp(i*pi)+0.0+0.0*i]; z2 = [z2 s2*exp(i*3/2*pi)-0.5+0.0*i].';
p2 = polygon([-0.5-0.5*i 0.0-0.5*i 0.0+0.0*i -0.5+0.0*i]);

z = [z1; z2];
u = -log(abs(z-c));                                     % Dirichlet boundary condition

N = 10+ceil(log(length(z)));                            % smooth part, nr. of terms
[a,H] = polyfitA(z-c,u,N);                              % Vandermonde with Arnoldi
wz = polyvalA(a,H,z-c);                                 % smooth part

up = u-real(wz);                                        % singular part of b.c.
[r,pol] = aaa(up,z,'mmax',1000,'lawson',0);             % AAA approximation
m = find((isinpoly(pol,p1,1E-16)==0) | ...              % find poles outside
(isinpoly(pol,p2,1E-16)==1));
B = 1./(z-pol(m).');                                    % coefficient matrix
B(:,end+1) = 1;
B(:,end+1) = 0; B(length(z1)+1:end,end) = 1;            % add column for modulus
V = [real(B) -imag(B)];
b = reshape(V\up,[],2);                                 % LS approximation
b = b(:,1)+i*b(:,2);
modulus = exp(-b(end))
maxerr = norm(u-real(B*b+wz),'inf')                     % maximum error at the boundary
b(end) = 0;
wp = B*b;                                               % singular part

w = wp+wz;                                              % solution

g = (z-c).*exp(w);                                      % conformal mapping function
f = aaa(g,z);                                           % polygon onto disk, AAA approx.
finv = aaa(z,g);                                        % disk onto polygon, AAA approx.
\end{lstlisting}
\begin{figure}[!h]
	\centering
	\includegraphics[width=0.75\linewidth]{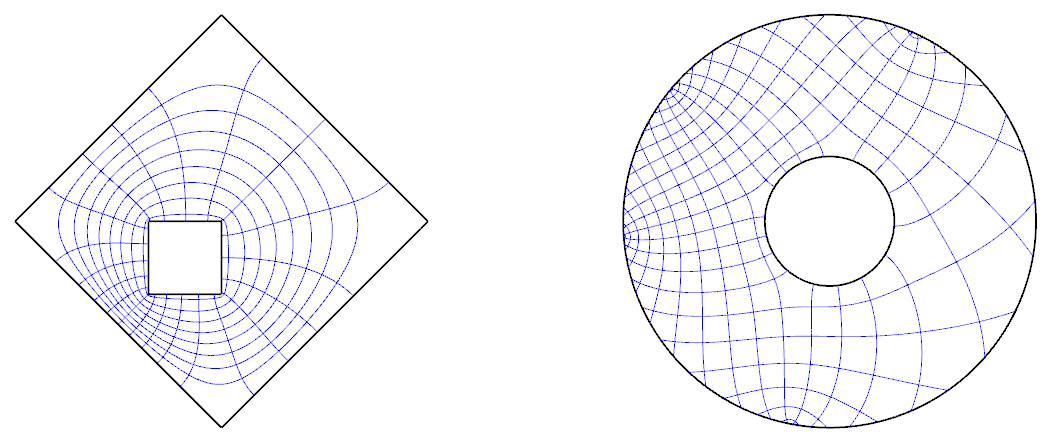}
	%\caption{}
	\label{fig:cm3_ann}
\end{figure}
Conformal mapping (see \cite{Tr19}) with uniform discretization step of 0.01 at both Jordan curves, from annulus of modulus $\rho=$0.314 onto polygon with hole (left), and from polygon onto annulus (the image lines of a square grid, right). $N$ = 17 in the smooth part and 246 poles, of which 108 fall outside the polygon, i.e. outside the outer square or inside the inner one. $|u-\Re w|_\infty$ = 4.78e-5.

\newpage
\section{Discussion}
The method described is conceived to be general-purpose, and clearly can not compete in speed nor accuracy with the specialized lightning Laplace solver, even so all the problems addressed, here and in other unreported tests, are solved in a few seconds on a reasonably modern laptop PC with good accuracy, nearly all of the time being consumed by the AAA approximation process. We emphasize again the fact that none of the solutions presented makes use of clustered samples in the neighbourhood of corners: while this clearly remains a possibility for the user, the AAA algorithm has always proven to be so flexible as to return decent approximations for any reasonably discretized contour. The same happens in cases of functions, not presented in this paper, that have singularities on smooth boundaries: the AAA algorithm clusters poles close to critical points following a similar pattern.

Besides the Laplace solver presented in this paper, and its possible applications to engineering computations, the idea of approximating analytic functions on adjacent regions by means of set of poles exterior to each other could suggest novel approaches in other topics: Riemann-Hilbert problems is one, where sectionally analytic functions are to be determined from boundary conditions.

\section*{Acknowledgements}
I am grateful to Lloyd N. Trefethen for his valuable suggestions that helped to improve this paper.

%\section*{References}
\bibliography{biblio.bib}
\bibliographystyle{elsarticle-num}
%\printbibliography

\end{document}